\documentclass{article} 

\usepackage{latexsym}
\usepackage[all]{xy}

\def\separation{\medskip}

\def\P{I\!\!P}

\def\O{{\cal O}}

\def\H{{\cal H}}

\def\I{{\cal I}}

\def\U{{\cal U}}
\def\V{{\cal V}}
\def\?{{\bf ??}}

\def\Hilb{{\rm Hilb}}

\def\dim{{\rm dim}}

\def\ker{{\rm ker}}

 \newtheorem{theorem}{Theorem}[section]
\newtheorem{lemma}[theorem]{Lemma} 
\newtheorem{prop}[theorem]{Proposition} 
 
\newtheorem{corollary}[theorem]{Corollary}
\newtheorem{remarks}[theorem]{Remarks} 
\newtheorem{remark}[theorem]{Remark}

\newcommand{\proof}{{\it Proof.}\ } 
\newcommand{\qed}{\hfill  $\Box$\separation}

\begin{document}

\title{The curve of lines on a prime Fano threefold \\
of genus 8}
\author{F. FLAMINI - E. SERNESI\footnote{Both authors are members of GNSAGA-INDAM.}}
\date{}
\maketitle

\section*{Introduction}

A complex projective three-dimensional nonsingular variety $X$ is a \emph{prime Fano threefold} if it 
has second Betti number $B_2=1$ and Pic$(X)$  is generated by $-K_X$. The (even) integer $(-K_X)^3$ is 
called the \emph{degree} of $X$ and 
\[
g := {1 \over 2}(-K_X)^3+1 \ge 2
\]
is the \emph{genus} of $X$.   
It is well known that prime Fano threefolds of genus $g$ exist only for   
$2 \le g \le 10$ or $g=12$.  Some of these threefolds   have generically non-injective period mapping  or even trivial intermediate jacobian despite the fact that they have non-trivial moduli. This is illustrated in \cite{sM04} for the case
$g=12$ and in \cite{IM07} for $g=8$.  On the other hand  a general prime Fano threefold $X$ of any genus has a Fano scheme of lines which is a nonsingular curve  
$\Gamma$ whose genus is known (\cite{enc99} Theorem 4.2.7).   Mukai has shown (see loc. cit. and references therein) that for $g=12$ a general  $X$ can be reconstructed from the  pair  $(\Gamma,\theta)$, where $\Gamma$ is the Fano curve (which has genus 3, being a plane quartic) and $\theta$ is a naturally defined  even theta-characteristic on $\Gamma$.  This is a modified version of the Torelli theorem where  the intermediate jacobian is replaced by   $(\Gamma,\theta)$.

It is natural to ask to which extent  this result   can be extended to other values of $g$, in the hope to have a good replacement of the intermediate jacobians by the 
jacobians of the Fano curves plus additional data, at least for Torelli-type purposes. 
This investigation seems to be still lacking in general.  

In the present paper we study the   genus $g=8$ case.  In this case the curve $\Gamma$ of lines contained in a 
general prime Fano threefold $X \subset \P^9$ has genus $26$. The reach geometry of such threefolds, 
investigated in \cite{gF30}, \cite{vI80}, \cite{pjP82}, \cite{nG93}, \cite{IM99}, \cite{sM02}, 
\cite{aK04} and \cite{IM07}, reflects into the geometry of $\Gamma$, which has a naturally defined even theta-characteristic $L$ which embeds it   in $\P^5$.   We prove the following result (Theorem \ref{T:main2}) analogous to the one proved by Mukai in the $g=12$ case:

\begin{theorem}
A general prime Fano threefold  $X$ of genus 8 can be reconstructed from the pair $(\Gamma,L)$, where $\Gamma$ is its Fano curve of lines and $L=\O_\Gamma(1)$ is the theta-characteristic   which gives the natural embedding $\Gamma \subset \P^5$. 
\end{theorem} 

In outline our proof goes as follows.  Every Fano threefold  $X$ of genus 8 can be realized as a linear section of  the grassmannian of lines of $\P^5$, in its Pl\"ucker embedding in $\P^{14}$, with a $\P^9$. As such, $X$ parametrizes a three-dimensional family of lines of $\P^5$ whose union is a quartic hypersurface $W$, called the \emph{Palatini quartic} of $X$. This quartic  has a singular curve $\Gamma(W)$ which turns out to be isomorphic to $\Gamma$, the curve of lines of $X$. This gives the natural embedding of $\Gamma$ in $\P^5$. Our strategy consists in proving the theorem in two steps:   we prove first that $\Gamma(W) \subset \P^5$ uniquely determines $W$, and then we show that there is only one Fano threefold of genus 8 in $\P^9$ whose associated Palatini quartic is $W$.

\noindent
For the first step we start   by representing $W$ as the degeneracy locus of a map $\phi$ of rank-five vector bundles, and $\Gamma(W)$ as the  scheme $D_3(\phi)$ where its rank drops twice. From this representation we are able to construct a locally free resolution of the ideal sheaf 
$\I:=\I_{\Gamma(W)}\subset \O_{\P^5}$ by means of the Gulliksen-Negard complex (see \cite{GN72} and \cite{jW03}).  
This resolution allows us to compute the dimensions of all the homogeneous pieces of the ideal of 
$\Gamma(W)$.  In particular we discover that $H^0(\I(3))$ has dimension six, and that it is generated 
by the partial derivatives of (an equation of) $W$. Then we are reduced to reconstructing $W$ from its jacobian ideal. 
We are able to do it in a special degenerate example and then we introduce some deformation theory which 
does it generically. In the Appendix  we prove the same result, using another approach. 
Namely, we show that if a Palatini quartic $W$ has linearly independent 
second partial derivatives (i.e. if it is not apolar to any quadric) then it can be reconstructed 
from $\Gamma(W)$. Since the linear independence of second partials is an open property, it suffices to find 
one such $W$.  Then we reproduce the script of a Macaulay2 program which determines a Palatini quartic with 
random coefficients having linearly independent   second partials.

\noindent
The second step is obtained by studying the Fano scheme $F(W)$ of lines on $W$.  One of its irreducible components is $X$ itself but,    being of degree 736 when
  Pl\"ucker embedded, it might a priori contain other similar components.  We exclude this possibility  by explicitly describing the   24 lines  of $W$ passing through a general point $p \in W$ and showing that only one of them can belong to a Fano  threefold of genus $8$ defining $W$.  An ingredient of our analysis is the \emph{dual cubic threefold}  $Y$ of $X$, whose properties and relation with $X$ are well known.

The paper consists of four sections and of an Appendix.  In  \S 1 we study the locally free resolution of the singular curve of a Palatini quartic. 
In \S 2 we show that a certain degenerate Palatini quartic can be reconstructed from the vector space generated by its partial derivatives. From this fact, by studying  the deformation theory of Palatini quartics, we prove that a general $W$ can be reconstructed from $\Gamma(W)$. 
In \S 3 we study the lines contained in the Palatini quartic and we prove the uniqueness of $X$, given $W$. We also give a complete description of all the irreducible components of $F(W)$. 
In the final  \S 4 we put together all the above to conclude.   
\separation

\emph{Acknowledgements.}  We are very grateful to G. Ottaviani for calling our attention on his paper 
\cite{gO92} and for providing us with the Macaulay program reproduced in our Appendix. We also thank 
E. Mezzetti,   M.L. Fania, D. Faenzi and E. Arrondo for valuable bibliographical references.  We thank  the referee for remarks  which contributed to improve the paper.

\section{The geometry of $\Gamma$}\label{S:gamma}

Let $V$ be a 6-dimensional complex vector space. Consider a  general Fano 3-fold of genus 8 and index 1
\[
X= G(2,V) \cap \P^9 \subset \P(\wedge^2V) \cong  \P^{14} 
\]
complete intersection of the grassmannian of lines of $\P(V)$, Pl\"ucker-embedded in $\P(\wedge^2V)$, with a general $\P^9$.   
We know (see  \cite{enc99} and \cite{dM83})  that $X$ has degree 14 and sectional genus 8.  We also know that the Hilbert scheme of lines on $X$ is a nonsingular irreducible curve 
$\Gamma \subset G(2,\wedge^2V)$ and that 
\[
g(\Gamma) = 26
\]
The union of the lines contained in $X$ is a ruled surface
\[
R_X = \bigcup_{\ell \in \Gamma} \ell \subset X
\]
 Since $X \subset G(2,V)$, each $x \in X$ parametrizes a line in $\P(V)$ which we denote by   $\ell_x$.  The union 
 \[
 W:= \bigcup_{x \in X} \ell_x  \subset \P(V)
 \]
 is a quartic hypersurface    called the \emph{Palatini quartic}
 of $X$ (see  \cite{IM99} and \cite{IM07}). It has a singular curve $\Gamma(W)$ which turns out to be isomorphic 
to $\Gamma$. 

It is known that $\Gamma(W)$ has degree $25$ and that $\O_{\Gamma(W)}(1)$ is a theta-characteristic, i.e. 
$\omega_\Gamma \cong \O_\Gamma(2)$ (cf. e.g. \cite{AR96}, p. 177, or \cite{aK04}, Remark 2.13). 
These facts can be proved as follows.

\noindent 
  $X$ is the intersection of $G(2,V)$ with a 
codimension-five linear subspace, which is defined by a $5$-dimensional subspace 
$U \subset \bigwedge^2V^\vee = H^0(\P(V), \Omega^1(2))$.  This subspace defines a map of rank-5 vector bundles on $\P(V)$:
 \[
  \xymatrix{
 U\otimes\O_{\P(V)}  \ar[r]^-\phi &  \Omega^1_{\P(V)}(2) 
  }
  \]
 As usual we will denote  by $D_k(\phi)$ the closed subscheme of $\P(V)$ defined by the condition 
 rk$(\phi) \le k$. Then   $D_4(\phi)$ is  the Palatini quartic $W \subset \P(V)$. 
 The curve $\Gamma(W)$ is  $D_3(\phi)$.  Let's compute the classes of $W$ and $\Gamma(W)$ using Porteous formula. 
 
Using the Euler sequence we obtain:
\[
 c_t(\Omega^1_{\P(V)}) = (1-t)^6 = 1-6t+15t^2-20t^3+15t^4-6t^5
\]
Therefore (see \cite{wF84}, p. 55):
\[
\begin{array}{lll}
c_t(\Omega^1_{\P(V)}(2))=&
  (1+2t)^5 - 6(1+2t)^4 t + 15(1+2t)^3 t^2 -20(1+2t)^2 t^3 \\ \\
&+15(1+2t)t^4 - 6t^5  \\ \\
&= 1 + 4t + 7t^2 + 6 t^3 + \cdots
\end{array}
\]
Then we get:
\[
\deg(W) = c_1= 4
\]
and
\[
\deg(\Gamma) = \left|\matrix{c_2&c_3 \cr c_1&c_2}\right| = c_2^2-c_1c_3 = 25
\]
Restricting $\phi$ to $\Gamma(W)$ we obtain:
\begin{equation}\label{E:resol0}
\xymatrix{
0 \ar[r]& E \ar[r]& U\otimes\O_{\Gamma(W)} \ar[r]^-{\phi_\Gamma} & 
\Omega^1_{\P(V)|\Gamma(W)}(2)\ar[r]&F\ar[r]&0 
}
\end{equation}
where $E$ and $F$ are rank-two vector bundles on $\Gamma$.  Thus, from \cite{GG73}, p. 145, we have:
\[
\omega_{\Gamma(W)} = \omega_{\P(V)|\Gamma(W)}\otimes \det(F \otimes E^\vee) 
= \O_{\Gamma(W)}(2)
\]
This means that $ \O_{\Gamma(W)}(1)$  is a theta-characteristic, in particular $\Gamma(W)$  has genus 
\[
g(\Gamma(W)) =26
\]

 \begin{prop}\label{P:gamma1}
 The ideal sheaf $\I=\I_{\Gamma(W)} \subset \O_{\P(V)}$ has the following locally free resolution:
 {\tiny\begin{equation}\label{E:resol1}
 \xymatrix{
 0\ar[r] & \O(-8)\ar[r] & (\Omega^4)^{\oplus 5}\ar[r] & L_{4,1}(T(-2))\oplus\O(-4)^{\oplus 24}
\ar[r]&\Omega^1(-2)^{\oplus 5}\ar[r]& \I\ar[r]&0 \\
  &&   U^{\vee} \otimes \Omega^4 \ar@{=}[u] && U^{\vee} \otimes \Omega^1(-2)\ar@{=}[u]}
 \end{equation}}where  $\O(-4)^{\oplus 24} = \frac{U^{\vee} \otimes U }{\bigwedge^5 U} \otimes \O(-4) $ and  
$L_{4,1}(T(-2))$ denotes the locally free sheaf of rank $24$ obtained by applying  the Schur functor $L_{4,1}$   to 
 $T(-2) =T_{\P(V)}(-2)= [ \Omega^1_{\P(V)}(2)]^\vee$, and which fits into an exact sequence as follows:
  \begin{equation}\label{E:resol2}
 \xymatrix{
 0 \ar[r] &\O(-4)\ar[r]\ar@{=}[d]&\Omega^4(4)\otimes \Omega^1(-2) \ar@{=}[d]\ar[r]&
 L_{4,1}(T(-2))\ar[r]&0\\
 &\omega^{-1}(-10)&T(-2)\otimes\bigwedge^4[T(-2)]}
 \end{equation}
 \end{prop}

 \proof (\ref{E:resol1}) is the Gulliksen-Negard 
resolution (see \cite{GN72}, \cite{jW03} n. 6.1.8), 
and the exact sequence (\ref{E:resol2}) makes explicit the fact that  $L_{4,1}(T(-2))$ is obtained by applying 
to $T(-2)$ the Schur functor relative to the partition   $(4,1)$.  \qed
 
 \begin{prop}\label{P:gamma2}
 $\Gamma(W)$ is linearly normal, and $k$-normal for all $k \ge 3$. It is not contained in any quadric and has $2$-deficiency  $h^1(\I_{\Gamma(W)}(2)) = 5$.  
  \end{prop}

  \proof 
  Tensoring by $\Omega^4(4)$  the twisted and dualized Euler sequence:
  \[
  \xymatrix{
  0\ar[r]& \Omega^1(-2) \ar[r]& V^{\vee} \otimes \O(-3)\ar[r] & \O(-2) \ar[r] & 0 }
  \]
  and recalling (\ref{E:resol2}) we obtain:
 \begin{equation}\label{E:resol7}
  \xymatrix{
  &0\ar[d] \\  &\O(-4)\ar[d] \\ 
  0\ar[r]&\Omega^4(4)\otimes \Omega^1(-2)\ar[d] \ar[r]& V^{\vee} \otimes \Omega^4(1)\ar[r] & \Omega^4(2) \ar[r] & 0 
  \\ &L_{4,1}(T(-2))\ar[d] \\ &0}
  \end{equation}
The exact row and Bott's theorem imply that 
$H^2(\Omega^4(4)\otimes \Omega^1(-2)(k))=0$ for all $k$.  Therefore the exact column gives 

\begin{equation}\label{E:resol3}
H^2(L_{4,1}(T(-2))(k)) =0  \qquad\hbox{ for all $k$}
 \end{equation}
  Now decompose the resolution (\ref{E:resol1}) into short exact sequences as follows:
  \[
  \xymatrix{
  0\ar[r] & K_1 \ar[r] & U^{\vee} \otimes \Omega^1(-2) \ar[r]& \I\ar[r]&0}
   \]
   \[
   \xymatrix{
  0 \ar[r]& K_2 \ar[r]&L_{4,1}(T(-2)) \oplus\O(-4)^{\oplus 24}\ar[r]& K_1 \ar[r] &0}
  \]
  \[
  \xymatrix{
 0\ar[r]&\O(-8)\ar[r]& U^{\vee} \otimes \Omega^4 \ar[r] & K_2 \ar[r] &0 }
  \]
   Using (\ref{E:resol3}) and chasing these exact sequences we see that  $H^2(K_1(k)) =0$ for all $k$.  
It follows that we have 
  a surjection:
  \[
  \xymatrix{
  U^{\vee} \otimes H^1(\Omega^1(k-2))  \ar[r]& H^1(\I(k)) \ar[r] & 0}
  \]
 for all $k$.  Using Bott's theorem again we deduce that  $H^1(\I(k)) = 0$ for all $k \ne 2$, and that we have a surjection:
 \[
 \xymatrix{
 {\bf C}^5 = U^{\vee} \otimes H^1(\Omega^1) \ar[r]& H^1(\I(2)) \ar[r] & 0}
 \]
 and therefore $h^1(\I(2)) \le 5$.  But, on the other hand 
 \[
 h^1(\I(2)) = h^0(\Gamma(W), \O_{\Gamma(W)}(2)) - h^0(\P(V),\O(2)) + h^0(\P(V),\I(2)) \ge 26-21=5
 \]
 Therefore $h^1(\I(2))=5$.   \qed

 \begin{corollary}\label{C:gamma1}
 The vector space 
  \[
  I_3 :=H^0(\I(3)) \subset H^0(\O_{\P(V)}(3))
  \]
has dimension 6 and coincides with the \emph{jacobian space} of $W$:   
\[
JW:=\left\langle {\partial W \over \partial X_0}, \dots, {\partial W \over \partial X_5}\right\rangle
\] 
  generated by the partial derivatives of   $W$ with respect to any system of homogeneous coordinated $X_0,\dots,X_5$ in $\P(V)$.  
 \end{corollary}

 \proof
 By Proposition \ref{P:gamma2} we have
  \[
 h^0(\I(3)) =  h^0(\P(V),\O(3)) - h^0(\Gamma,\O_\Gamma(3)) =  6
 \]
 Since $W$ is singular along  $\Gamma(W)$, it follows that $I_3$ contains the space generated by the six partial derivatives of an equation of $W$.  Since  $W$ is not a cone it follows that  the six partial derivatives of $W$ are linearly independent. \qed

 \section{Reconstruction of $W$ from $\Gamma(W)$}\label{S:family}
 
  The homomorphisms $\phi: U\otimes\O_{\P(V)} \to \Omega^1_{\P(V)}(2)$ defining the Palatini quartics
 are parametrized by an open set of the grassmannian $G(5,\bigwedge^2V^\vee)$, which is   irreducible and nonsingular of dimension 50. The rule
 \[
 \xymatrix{
 \phi\ar@{|->}[r]& W=D_4(\phi)}
 \]
 defines a rational map
 \[
  \xymatrix{w:G(5,\bigwedge^2V^\vee)\ar@{-->}[r]& \P(S^4V^\vee)}
  \]
 whose image is the locally closed subset     $\U \subset \P(S^4V^\vee)$  parametrizing Palatini quartics. It follows that $\U$ is irreducible of dimension $\le 50$.

 Consider a Palatini quartic $W\subset \P(V)$.  From Proposition \ref{P:gamma2} we deduce that
 \[
 h^0(\I(4)) = h^0(\O_{\P(V)}(4)) - h^0(\O_\Gamma(4)) = 126 - 75 = 51. 
 \]Observe that we have an exact sequence:
 \begin{equation}\label{E:addrem0}
 \xymatrix{
 0\ar[r] & \O_{\P(V)}\ar[r] & \I(4) \ar[r] & N'_W \ar[r] & 0}
 \end{equation}  
  where $N'_W \subset N_W=\O_W(4)$ is 
 the {\em equisingular normal sheaf} (\cite{eS06}, Proposition 1.1.9).   
 From (\ref{E:addrem0}) we obtain:
 \[
 h^0(N'_W) = 50,  \qquad h^1(N'_W)=0
 \]
     Therefore from local deformation theory we deduce that  the locally trivial deformations of $W$ are unobstructed and of dimension 50 (\cite{eS06}, Example 4.7.1(i)).

   \begin{lemma}\label{L:addrem0}  Assume that
   \[
  \xymatrix{
 U\otimes  \O_{\P(V)}  \ar[r]^-\phi &  \Omega^1_{\P(V)}(2)
  }
  \]
  is a homomorphism defining a Palatini quartic $W \subset \P(V)$. Then 
   \[
   dw_\phi: \xymatrix{
  T_\phi G(5,\bigwedge^2V^\vee)\ar[r] & H^0(N_W)}
  \]
    is injective and {\em Im}$(dw_\phi)= H^0(N'_W)$.    Therefore $w$    is unramified at $\phi \in G(5,\bigwedge^2V^\vee)$ and  $\U$, with its reduced scheme structure,  is irreducible  and nonsingular of dimension $50$.   
  \end{lemma}
  
  \proof  
      We have:
    \[
    T_\phi G(5,\bigwedge^2V^\vee) = U^\vee\otimes\left[ (\bigwedge^2V^\vee)/U\right] =
    [U^\vee\otimes H^0(\Omega^1(2))]/[U^\vee\otimes U]
    \]
    and $dw_\phi$ is induced by   the composition:
    \[
    \xymatrix{
    U^\vee\otimes H^0(\Omega^1(2))\ar[r]& H^0(\I(4))\ar[r]&  H^0(N'_W)}
    \]
    where the first map  comes from   the last map to the right in the resolution (\ref{E:resol1}) twisted by $\O(4)$. Therefore  Im$(dw_\phi) \subset H^0(N'_W)$ and to prove equality it suffices to show that 
     \[
    \xymatrix{
    U^\vee\otimes H^0(\Omega^1(2))\ar[r]& H^0(\I(4))}
    \]
    is surjective. 
    Decompose the resolution (\ref{E:resol1}) into short exact sequences as in the proof of Proposition \ref{P:gamma2} and twist everything by $\O(4)$.  We obtain the exact sequences:
  \begin{equation}\label{E:resol4}
  \xymatrix{
  0\ar[r] & K_1(4) \ar[r] & U^{\vee} \otimes \Omega^1(2) \ar[r]& \I(4)\ar[r]&0}
  \end{equation}
  \begin{equation}\label{E:resol5}
  \xymatrix{
  0 \ar[r]& K_2(4) \ar[r]&[L_{4,1}(T(-2))](4) \oplus\O^{\oplus 24}\ar[r]& K_1(4) \ar[r] &0}
  \end{equation}
  \begin{equation}\label{E:resol6}
  \xymatrix{
  0\ar[r]&\O(-4)\ar[r]& U^{\vee} \otimes \Omega^4(4) \ar[r] & K_2(4) \ar[r] &0 }
  \end{equation}
  From (\ref{E:resol4}) we see that  it suffices to show that 
  \[
  H^1(K_1(4))=0 
  \]
  From (\ref{E:resol6}) and Bott's Theorem we obtain $H^i(K_2(4))=0$ for $i=0,1,2$.  From (\ref{E:resol5})
  we deduce:
  \[
  h^1(K_1(4)) = h^1(L_{4,1}(T(-2))(4))
  \]
  so that we are reduced to show that:
  \begin{equation}\label{E:resol8}
   H^1(L_{4,1}(T(-2))(4)) = 0
  \end{equation}
  We consider the diagram (\ref{E:resol7}) twisted by $\O(4)$:
  \begin{equation}\label{E:resol9}
  \xymatrix{
  &0\ar[d] \\  &\O\ar[d] \\ 
  0\ar[r]&\Omega^4(8)\otimes \Omega^1(2)\ar[d] \ar[r]& V \otimes\Omega^4(5) \ar[r] & \Omega^4(6) \ar[r] & 0 
  \\ &L_{4,1}(T(-2))(4)\ar[d] \\ &0}
  \end{equation}
  Bott's Theorem gives $H^i(\Omega^4(5-i))=0$ for all $i \ge 1$, so that $\Omega^4$ is 5-regular. Therefore the map
  \[
  \xymatrix{ V \otimes H^0(\Omega^4(5)) \ar[r] & H^0(\Omega^4(6))}
  \]
  is surjective.  Since, again by Bott's Theorem:
  \[
  H^1(\Omega^4(5))= 0
  \]
  we obtain
  \[
  H^1(\Omega^4(8)\otimes \Omega^1(2)) = 0
  \]
  Finally, from the column of  (\ref{E:resol9}) we obtain  (\ref{E:resol8}) and the proof of 
  Im$(dw_\phi)= H^0(N'_W)$ is completed.
  
  \noindent
  What this says is that  $\U$ is supported, locally at each $W$,    on the locus of locally trivial deformations of $\U$, which is nonsingular of dimension $50$, as observed above.    \qed

  The rules:
 \[
 \xymatrix{
 \phi\ar@{|->}[r]\ar@{|->}[d]& W=D_4(\phi)\ar@{|->}[dl] \\
 \Gamma=D_3(\phi)={\rm Sing}(W)}
 \]
 define rational maps which fit into a commutative diagram:
 \[
 \xymatrix{G(5,\bigwedge^2V^\vee)\ar@{-->}[r]^-w\ar@{-->}[d]_-\gamma & \U\ar[dl]^-\sigma \\
\H}
 \]
  where $\H$ is an open subscheme of $\Hilb^{\P(V)}_{25t-25}$ containing the points parametrizing  the curves Sing$(W)$. It is easy to compute that dim$(\H) \ge 100$ (see Remark \ref{R:main1}). This means that the points parametrizing the curves ${\rm Sing}(W)$, as $W \in \U$,  fill a locally closed subscheme of   codimension $\ge 50$ in $\Hilb^{\P(V)}_{25t-25}$.
 The singular curve $\Gamma(W) = {\rm Sing}(W)$ of a Palatini quartic $W$ is 3-normal and satisfies 
 $h^0(\I(3))=6$ (Proposition \ref{P:gamma2}) and this property is open in $\H$  by semicontinuity.  Therefore we have a well defined  rational map
 \[
 i_3:\xymatrix{
 \H \ar@{-->}[r] & G(6,S^3V^\vee) }
  \]
  defined by $i_3(C) = H^0(\I_C(3))\subset S^3V^\vee$.   This map fits into the previous diagram as follows:
  \begin{equation}\label{E:diagram1}
   \xymatrix{G(5,\bigwedge^2V^\vee)\ar@{-->}[r]^-w\ar@{-->}[d]_-\gamma & \U\ar[dl]^-\sigma\ar[d]^-J \\
 \H \ar@{-->}[r]^-{i_3} & G(6,S^3V^\vee)}
 \end{equation}
 where $J$ is the morphism associating to a Palatini quartic its jacobian vector space. The commutativity relation $J = i_3\circ \sigma$ is proved in Corollary \ref{C:gamma1}.

  From the commutative diagram  (\ref{E:diagram1}) we deduce immediately: 
  
  \begin{lemma}\label{L:addrem1}
    If  $J$ is generically injective  then $\sigma$ is generically injective.  
     \end{lemma}

    We will prove that $J$ is generically injective by a degeneration  argument. 
      Observe that  
         $J$ is the restriction to $\U$ of a morphism:
     \[
     \tilde J: \xymatrix{ \V\ar[r] & G(6,S^3V^\vee)}
 \]
  where $\V \subset \P(S^4V^\vee)$ is the open subset consisting of the hypersurfaces $F$  whose jacobian vector space has dimension six.

  \begin{lemma}\label{L:donagi}
  The morphism $\tilde J$ is birational onto its image and has connected fibres.
  \end{lemma}
  
  \proof
  Let   $F,G \in \V$ be such that $\tilde JF=\tilde JG$, and let  $F_t=tF+(1-t)G$,  $t \in {\bf C}$. Then for each $t \in {\bf C}$ and $j=0,\dots, 5$ we have:
\[
{\partial F_t \over \partial X_j} = t   {\partial F \over \partial X_j}+ (1-t){\partial G\over \partial X_j} \in
\tilde JF+\tilde JG= \tilde JF
\]
 Therefore   the jacobian vector space of   $F_t$ is contained in $\tilde JF$. It follows  that   there is an open subset $\emptyset \ne A \subset {\bf C}$ such that $F_t \in \V$ and $\tilde JF_t=\tilde JF$ for all $t \in A$. This proves that the fibres of $\tilde J$ are connected.
  
  \noindent
  By \cite{CG80},  Lemma p. 72, if $F,G \in \V$ are such that $\tilde JF=\tilde JG$ and $F$ is general, then $F=G$: this implies that  $\tilde J$ is birational onto its image.  \qed
  
  \emph{Note:} the proof of connectedness of the fibres  is an adaptation of an argument which appears in the proof of Proposition 1.1 of \cite{rD83}.

   Before we state the next proposition, observe   that a homomorphism 
 \[
  \xymatrix{
  U \otimes \O_{\P(V)}  \ar[r]^-\phi &  \Omega^1_{\P(V)}(2)
  }
  \]
  is defined by assigning five linearly independent elements of 
  \[
  H^0(\Omega^1_{\P(V)}(2))= \bigwedge^2V^\vee
  \]
 i.e., after choosing a basis of $V$, by  five linearly independent $6 \times 6$ skew-symmetric matrices.

  \begin{prop}\label{P:matrix}
 Consider the following linearly independent skew-symmetric matrices:
 \[
 M_1=\pmatrix{0&1&0&0&0&0 \cr -1&0&0&0&0&0 \cr 0&0&0&0&0&0 \cr
 0&0&0&0&0&0 \cr 0&0&0&0&0&0 \cr 0&0&0&0&0&0 }\qquad
 M_2=\pmatrix{0&0&0&0&0&-1\cr 0&0&1&0&0&0 \cr 0&-1&0&0&0&0\cr
 0&0&0&0&0&0 \cr 0&0&0&0&0&0 \cr 1&0&0&0&0&0 }
 \]
 \[
 M_3=\pmatrix{0&0&0&0&0&0 \cr 0&0&0&0&0&0 \cr  0&0&0&1&0&0 \cr
  0&0&-1&0&0&0 \cr 0&0&0&0&0&0  \cr 0&0&0&0&0&0} \qquad
  M_4 =\pmatrix{0&0&0&0&0&0 \cr 0&0&0&0&0&0 \cr 0&0&0&0&0&0 \cr
  0&0&0&0&1&0 \cr 0&0&0&-1&0&0  \cr 0&0&0&0&0&0}
  \]
  \[
  M_5 = \pmatrix{0&0&0&0&0&0 \cr
 0&0&0&0&0&0 \cr 0&0&0&0&0&0 \cr 0&0&0&0&0&0 \cr 0&0&0&0&0&1 \cr0&0&0&0&-1&0}
  \]
 Then:

 \begin{description}
 \item[(i)]
 The degeneracy hypersurface   of the homomorphism
 \[
  \xymatrix{
  \phi_M:  \O_{\P^5}^{\oplus 5}  \ar[r] &  \Omega^1_{\P^5}(2)
  }
  \]
 defined by the 5-tuple $M = (M_1, \dots, M_5)$  has equation:
 \[
 W_M: \quad X_1X_2X_3X_4 - X_0X_3X_4X_5 = 0
 \]
  \item[(ii)] The partial derivatives of $W_M$ are linearly independent, and $W_M$ is the only quartic 
hypersurface whose partial derivatives generate the space
  \[
 \left\langle {\partial W_M \over \partial X_0}, \dots, {\partial W_M \over \partial X_5}\right\rangle
\] 
 \end{description}
\end{prop}
  
 \proof   (i) Clearly   $D_4(\phi_M)$ is defined by the six  $5\times 5$ minors of the matrix:
 \[
 A= \pmatrix{
 X_1&-X_0&0&0&0&0 \cr -X_5&X_2&-X_1&0&0&X_0 \cr 0&0&X_3&-X_2&0&0 \cr
 0&0&0& X_4&-X_3&0 \cr 0&0&0&0&X_5&-X_4}
 \]
  and an easy computation shows that the $i$-th such minor is:
  \[
  X_i( X_1X_2X_3X_4 - X_0X_3X_4X_5), \qquad  i=0,\dots, 5
  \]
   This proves (i).
   
   (ii) The partial derivatives of $W_M$ are:
   \[
   {\partial W_M \over \partial X_0} = X_3X_4X_5, \qquad  
   {\partial W_M \over \partial X_1} = X_2X_3X_4, \qquad 
    {\partial W_M \over \partial X_2} = X_1X_3X_4
   \]
   \[
    {\partial W_M \over \partial X_3} = X_4(X_1X_2-X_0X_5), \qquad 
     {\partial W_M \over \partial X_4} = X_3(X_1X_2-X_0X_5)
     \]
     \[ 
     {\partial W_M \over \partial X_5} = -X_0X_3X_4
     \]
       Their independence is obvious.  In order to prove the uniqueness of $W_M$  
        it suffices to show that if
   $\{F_0,\dots,F_5\}$ is a basis of the space  
   $ \left\langle {\partial W_M \over \partial X_0}, \dots, {\partial W_M \over \partial X_5}\right\rangle$ satisfying
   \begin{equation}\label{E:partials}
   {\partial F_i \over \partial X_j} =  {\partial F_j \over \partial X_i}, \qquad 0 \le i,j \le 5 
   \end{equation}
   then   
   \[
   F_i = \lambda {\partial W_M \over \partial X_i}, \qquad i=0,\dots, 5
   \]
    for some $\lambda \ne 0$.  We follow the method introduced in \cite{CG80}  (Lemma p. 72), which consists in writing 
    \[
    F_i = \sum_\alpha A_{i\alpha}  {\partial W_M \over \partial X_\alpha},  \qquad  i=0,\dots, 5
    \]
    and imposing the conditions (\ref{E:partials}). This leads to the identity:
    \begin{equation}\label{E:hessian}
    AH = H \ {^t\!A}
    \end{equation}
    where 
    \[
    H =  \pmatrix{{\partial^2 W_M \over \partial X_i\partial X_j}}
    \]
    is the hessian matrix of $W_M$.  Now a direct computation shows that the only $6\times 6$ matrices $A$ satisfying (\ref{E:hessian}) are of the form  $A = \lambda \ {\bf I}_6$ for some constant 
    $\lambda \ne 0$.  \qed

    Observe that $W_M$ is not a Palatini quartic, because it is a reducible hypersurface, but it belongs to the closure of $\U$.

    \begin{corollary}\label{C:gentor2}  $J$ is generically injective and therefore $\sigma$ is generically injective. 
    \end{corollary}

    \proof   Let $\widetilde\V\subset \V$ be the open set where $\tilde J$ is injective.  Proposition \ref{P:matrix}    implies that $W_M \in \widetilde \V$.  Since $W_M$ belongs to the closure of $\U$, which is irreducible, by the upper-semicontinuity of the dimension of the fibres of $\widetilde J$ and by Lemma \ref{L:donagi}   we deduce that 
    $\U \cap \widetilde\V \ne \emptyset$.
This proves that $J$ is generically injective. 
The generic injectivity of $\sigma$ follows from Lemma \ref{L:addrem1}. \qed

\section{Lines on the Palatini quartic}\label{S:W}

If $Z \subset \P^n$ is a hypersurface, the  Hilbert scheme of lines contained in $Z$,  which is usually called the \emph{Fano scheme} of lines on $Z$, will be denoted  by $F(Z)$.

\begin{prop}\label{P:pala1}
Let $Z \subset \P^4$ be a quartic hypersurface with at most isolated singularities, not a cone  and not containing a $\P^2$. Then $F(Z)$ is purely one-dimensional.
\end{prop}

\proof  It is well known (\cite{wF84}, Example 14.7.13 p. 275) that $F(Z)$ is the vanishing scheme of a
section of $S^4Q$, where $Q$ is the universal quotient bundle on $G(2,5)$. Therefore each component of $F(Z)$ has dimension $\ge 1$.  Assume that $[\ell] \in F(Z)$ is a general point of a component $T$ of dimension $\ge 2$.   If $\ell\cap {\rm Sing}(Z)=\emptyset$ then we have an exact sequence
\[
\xymatrix{
0\ar[r]&N_{\ell/Z} \ar[r]& \O_\ell(1)^3 \ar[r]& \O_\ell(4)\ar[r] & 0 }
\]
which implies that $N_{\ell/Z} \cong \O_\ell(a)\oplus \O_\ell(b)$ with $a<0$ and $b>0$. 
It follows that $\ell$ is not free and therefore the union of the lines of $T$ is a surface $S$ contained in $Z$. But the only surface containing a two-dimensional family of lines is $\P^2$, therefore $S=\P^2$, and this contradicts 
the hypothesis. 

\noindent
If  $\ell\cap {\rm Sing}(Z)\not=\emptyset$  then all the lines of $T$ contain a fixed singular point of $Z$ and, since $Z$ is not a cone,  the union of the lines of $T$ must be a surface $S$ contained in $Z$. We conclude as before.
\qed

We keep the notation of \S \ref{S:gamma}.  In particular, 
\[
X=G(2,V)\cap \P(U)^\perp\subset \P(\wedge^2V)
\]is a general prime Fano threefold of genus 8, where $U \subset  \wedge^2V^\vee$ is a vector subspace of dimension five, 
and $W \subset \P(V)$ is the Palatini quartic associated to $X$. 
Consider the incidence relation
\[
\widetilde W := \{(x,p): p \in \ell_x \}  \subset X \times \P(V)
\]
Clearly  $\widetilde W $ is a $\P^1$-bundle over $X$, in particular it is a nonsingular fourfold:  in fact $\widetilde W = \P({\cal T})$ where  ${\cal T}$ is the restriction to $X$ of the tautological rank-2 vector bundle on $G(2,V)$.  The second projection maps $\widetilde W$  onto  $W$.  Denote by
\[
\xymatrix{
&\widetilde W \ar[dl]_{p_X} \ar[dr]^{p_W} \\
X && W}
\]
 the morphisms induced by the projections.

 \begin{prop}\label{P:lines1}
 \begin{description}
 
 \item[(i)]
 The restriction
 \[
 \xymatrix{ p_W : \widetilde W \backslash p_W^{-1}(\Gamma(W)) \ar[r] & W \backslash \Gamma(W)}
 \]
 is an isomorphism.
 
 \item[(ii)]  
 $\xymatrix{p_W^{-1}(\Gamma(W)) \ar[r]^-{p_W} & \Gamma(W)}$ is a ruled surface. 
 
 \item[(iii)] For each $p \in \Gamma(W)$
 \[
 \ell :=p_X(p_W^{-1}(p)) \subset X
 \] is a line.
 
 \item[(iv)]  Conversely, for each line $\ell \subset X$ we have
 \[
 p_X^{-1}(\ell) \cong F_1 := \P(\O\oplus \O(-1))
 \]
 and $p_W$ contracts the negative section of  $p_X^{-1}(\ell)$ to a point $p \in \Gamma(W)$, and the fibres of  $p_X^{-1}(\ell) \to \ell$  are mapped by $p_W$ to the lines of a pencil centered at $p$ and spanning a plane $\Pi_p$.
 \end{description}

 \end{prop}
 
 \proof  See    \cite{aK04}, \cite{IM07}.  \qed

  \begin{corollary}\label{C:lines1}
  \begin{description}
 
 \item[(i)]  $p_W: \widetilde W \to W$ is a small contraction,

 \item[(ii)]  $p_W(p_X^{-1}(R_X)) \subset W$ is an irreducible Weil divisor which is the union of the one-dimensional family of planes  $\Pi_p$,  $p \in \Gamma$.

 \item[(iii)] Each point $p \notin \Gamma(W)$ belongs to a unique line $\ell_x$, for some $x \in X$.
 
 \item[(iv)] For each $x \in X \backslash R_X$ we have 
 \[
 N_{\ell_x/W}  \cong \O_{\ell_x}^{\oplus 3}
 \]
  \end{description}
 \end{corollary}

 \proof (i),(ii) and (iii) are an immediate consequence of the proposition. (iv) follows from the fact that 
$\ell_x =p_W (p^{-1}_X(x)) \subset W$ and $p_X^{-1}(x)\subset \widetilde W$ have isomorphic open neighborhoods, 
and therefore isomorphic normal bundles. Since $p_X^{-1}(x)$ is a fibre of a fibration, it has trivial normal bundle.    \qed

   We now introduce the   \emph{dual cubic threefold} of $X$, which is defined  as follows. 
The 5-dimensional vector space $U \subset \bigwedge^2V^\vee = H^0(\P(V), \Omega^1(2))$ 
 can be also interpreted as defining a linear map:
\[
\xymatrix{
U\otimes V \ar[r] & V^\vee}
\]
which induces a homomorphism:
\[
\xymatrix{
\zeta: V\otimes \O_{\P(U)}(-1) \ar[r] & V^\vee \otimes \O_{\P(U)}}
\]
of   free sheaves of rank six on $\P(U) \cong \P^4$.  Because of the skew-symmetry of the elements of $U$ the degeneracy locus of $\zeta$ is a cubic threefold $Y \subset \P(U)$ defined by the pfaffian of $\zeta$, and we have an exact sequence:
\[
\xymatrix{
0 \ar[r] & V\otimes \O_{\P(U)}(-1) \ar[r]^-\zeta & V^\vee \otimes \O_{\P(U)}\ar[r] & E \ar[r] & 0}
\]
where $E$ is a  rank-two locally free sheaf on $Y$. We have morphisms:
\[
\xymatrix{& \P(E) \ar[dl]_{q_Y} \ar[dr]^{q_W} \\
Y&&W}
\]
 where $q_Y$ is the natural projection.   $q_W$ is defined by viewing  each    fibre $E_y$ of $E$, $y \in Y$,  as a 2-vector subspace $E^\vee_y \subset V$  thus defining a line in $\P(V)$, called a \emph{kernel line}, and denoted  $n_y$. It can be shown that each kernel line is contained in $W$. 
 
 \begin{prop}\label{P:lines4}
 \begin{description}
 
 \item[(i)] The restriction
 \[
 \xymatrix{ q_W : \P(E) \backslash q_W^{-1}(\Gamma(W)) \ar[r] & W \backslash \Gamma(W)}
 \]
 is an isomorphism.
 
 \item[(ii)]  
 $\xymatrix{q_W^{-1}(\Gamma(W)) \ar[r]^-{q_W} & \Gamma(W)}$ is a ruled surface. 
 
 \item[(iii)] For each $p \in \Gamma(W)$
 \[
 q_Y(q_W^{-1}(p)) \subset Y
 \] 
 is a line   such that   $E_{|q_Y(q_W^{-1}(p))} \cong \O\oplus\O(-2)$.

 \item[(iv)]  Conversely, for each line $\lambda \subset Y$ 
 such that $E_{|\lambda} \cong \O\oplus\O(-2)$ we have
 \[
 q_Y^{-1}(\lambda) \cong F_2 := \P(\O\oplus \O(-2))
 \]
 and $q_W$ contracts the negative section of  $q_Y^{-1}(\lambda)$ to a point $p \in \Gamma(W)$, and the fibres of  $q_Y^{-1}(\lambda) \to \lambda$  are mapped by $q_W$ to the lines of a quadric cone $K_p$ with vertex  $p$.
 
 \noindent
 The lines $\lambda \subset Y$ with this property are called \emph{B-lines}, and their union is a ruled surface $R_Y \subset Y$.   The other lines $\ell \subset Y$ are called  \emph{A-lines}.

 \end{description}
  \end{prop}
  
   \proof  See   \cite{AR96},  \cite{aK04}, \cite{IM07}.  \qed

  Analogously to the case of $\widetilde W$,  we have the  
  
 \begin{corollary}\label{C:lines2}
 \begin{description}
 
 \item[(i)]  $q_W: \P(E) \to W$ is a small contraction, and

 \item[(ii)]  $q_W(q_Y^{-1}(R_Y)) \subset W$ is an irreducible Weil divisor which is the union of the one-dimensional family of cones $K_p$,  $p \in \Gamma(W)$. 
 
 \item[(iii)]  Each point $p \notin \Gamma(W)$ belongs to a unique kernel line.
 
  \item[(iv)] Let $y \in Y \backslash q_Yp^{-1}_W(\Gamma(W))$.  Then
 \[
 N_{n_y/W}   \cong \O_{n_y}^{\oplus 3}
 \]
 \end{description}
 \end{corollary}
 
 \proof It is similar to the proof of Corollary \ref{C:lines1}.  \qed

  It is well known that $X$ and $Y$ are birationally equivalent (see \cite{enc99}, Theorem 4.6.6 and  \cite{kT89}).   From this fact we can deduce the following:
  
  \begin{lemma}\label{L:pala1}
  $W$ does not contain a $\P^3$.
  \end{lemma}
  
  \proof  Assume that there is a $\P^3$, say $\Lambda$,  contained in $W$. Then, since 
the general line $\ell_x$ is not contained in $\Lambda$, a general point of $\Lambda$ 
is contained in a unique line $\ell_x$ by Corollary \ref{C:lines1}(iii). Therefore we obtain a birational 
map $\xymatrix{X\ar@{-->}[r]& \Lambda = \P^3}$ so $X$ is rational, 
contradicting the fact that $X$ is birational to $Y$. \qed

We can now prove the following:

\begin{prop}\label{P:pala2}
  The Fano scheme  $F(W)$ is a local complete intersection of pure dimension three. In particular the general point of $W$ is contained in  finitely many lines contained in $W$.
\end{prop}

\proof
$F(W)$ is  the vanishing scheme of a section of $S^4Q$, where $Q$ is the universal quotient bundle on $G(2,V)$. Therefore every component of $F(W)$ has dimension $\ge 3$.  Since $W$ is not a cone and does not contain a $\P^3$    a general hyperplane section   $H \cap W$ is a quartic hypersurface satisfying the hypothesis of Prop. \ref{P:pala1}: therefore   its Fano scheme of lines $F(H \cap W)$ is purely one-dimensional. If $F(W)$ had a component $B$ of dimension $\ge 4$ then $B$ would give rise to a subscheme of $F(H\cap W)$ of dimension $\ge \dim(B)-2 \ge 2$ and this is a contradiction. Therefore  $F(W)$ has pure dimension three and, since it is the vanishing scheme of a section of a vector bundle of rank five, equal to its codimension, it follows that it is a local complete intersection. The last assertion is obvious.  \qed

Using the previous  results  we can prove the following:
  
  \begin{prop}\label{P:lines2}
  \begin{description}
  
 \item[(i)]  The classifying morphism $\xymatrix{\chi:X \ar[r] & F(W)}$ induced by the   family
  \begin{equation}\label{E:lines1}
  \xymatrix{
   \widetilde W \ar[d]_{p_X}\ \ar@{^(->}[r]&   X \times W \\
  X}
  \end{equation}
  is a closed embedding which, composed with the embedding $F(W)\subset G(2,V)$, gives the natural inclusion $X \subset G(2,V)$ and identifies $X$ with an irreducible component of $F(W)$. 
  
  \item[(ii)]
  The classifying morphism $\xymatrix{h:Y \ar[r] & F(W)}$ induced by the family
  \begin{equation}\label{E:lines7}
  \xymatrix{
   \P(E) \ar[d]_{q_Y}\ \ar@{^(->}[r]^-n&   Y \times W \\
  Y}
  \end{equation}
  where $n$ embeds every fibre of $q_Y$ as the corresponding kernel line,  is a closed embedding whose image is an irreducible  component of $F(W)$.  The composition 
  \[
  \xymatrix{
  Y \ar[r]^-h& F(W) \ \ar@{^(->}[r] & \P(\wedge^2V^\vee)}
  \]
  is given by the   linear system $|-2K_Y|$. In particular $h(Y)$ is   an irreducible component of $F(W)$ of degree 24.
  \end{description}
  \end{prop}

  \proof
  (i)
  Since we have a commutative diagram:
  \[
  \xymatrix{
  X \ar@{^(->}[rr]^-\iota\ar[dr]^-\chi && G(2,V) \\
  &F(W) \ar@{^(->}[ur]^-j}
  \]
where $\iota$ and $j$ are closed embeddings, it follows that $\chi$ is a closed embedding as well.  Moreover, by Corollary \ref{C:lines1}(iv)  the tangent space $T_xF(W)$ has dimension three at a general point $x \in X$ and this proves that $X$ is an irreducible component of $F(W)$. 

\noindent
(ii) As in the proof of Prop.  \ref{P:lines2} one shows that $h$ is a closed embedding.  The composition $\xymatrix{Y \ar[r]& F(W) \ \ar@{^(->}[r] & \P(\wedge^2V^\vee)}$ 
  is defined by the invertible sheaf $\wedge^2E$, which is isomorphic  to $\O(-2K_Y)$ (\cite{IM99}, Theorem 2.2).
  \qed

A closer analysis shows that $F(W)$ contains   other irreducible components (see Remark \ref{R:pala1}). The following result will be crucial.

\begin{prop}\label{P:pala3}
A general point of $W$ is contained in precisely $24$ distinct lines contained in $W$.  
\end{prop}

\proof  Let $p \in W$ be a general point. Then there are at most $4!=24$ distinct lines in $W$ containing $p$ \cite{jL03}.  Therefore it suffices to find at least 24 such lines.

\noindent
There is a unique $x \in X$ such that $p \in \ell_x$ and a unique $y \in Y$ such that $p \in n_y$, and $\ell_x \ne n_y$.  By the genericity of $p$ the points $x$ and $y$ are general. Therefore $x$ is contained in 16 distinct nonsingular conics $q_1,\dots, q_{16} \subset X$ \cite{kT89} and each of them is of type $(0,0)$, i.e. it has normal bundle isomorphic to
$\O\oplus \O$ (\cite{enc99}, Prop. 4.2.5).  To   such conics there correspond  nonsingular two-dimensional quadrics $Q_1, \dots , Q_{16} \subset W$ containing $\ell_x$, each of them being the union of the lines $\ell_{x'}$,  $x' \in q_i$.  For each $i=1, \dots, 16$ denote by 
$\lambda_i \subset Q_i$ the unique line containing $p$ and different from $\ell_x$.  The lines $\lambda_1, \dots, \lambda_{16}$ are distinct.

\noindent
Similarly, we have that $y$ is contained in six distinct  lines $l_1, \dots, l_6 \subset Y$ \cite{jL03} and they are A-lines. For each $j=1,\dots, 6$   we obtain a nonsingular quadric (\cite{IM07}, Lemma 2.2)
\[
Q^j = \bigcup_{y'\in l_j} n_{y'}
\]
Denote by $\nu_j \subset Q^j$ the unique line containing $p$ and different from $n_y$.  

\noindent
Altogether we have found 24 lines containing $p$ and contained in $W$:
\[
\begin{array}{ll}
\ell_x, n_y \\
\lambda_1, \dots, \lambda_{16} \\
\nu_1, \dots, \nu_6 
\end{array}
\]
Clearly $\ell_x, \lambda_1, \dots, \lambda_{16}$ are distinct and $n_y, \nu_1, \dots, \nu_6$ are distinct. 
Since $\ell_x \ne n_y$ it remains to be excluded that some $\lambda_i$ is a kernel line and that some $\nu_j$ is a 
$\ell_{x'}$, for some $x' \in X$.  But either possibility implies that $Q_i=Q^j$ for some $i,j$, and this is impossible because the quadrics $Q_1, \dots , Q_{16}, Q^1, \dots, Q^6$ are distinct (see \cite{IM07}, p. 511, for a proof of this fact). \qed

We can now prove the main result of this section.

\begin{theorem}\label{C:pala1}
$X$ is the only prime Fano threefold of genus 8 having $W$ as its associated Palatini quartic.
\end{theorem}

\proof We keep the notation of Prop. \ref{P:pala3}.  Let $p\in W$ be a general point, and let  
\[
\ell_x, n_y,
\lambda_1, \dots, \lambda_{16},
\nu_1, \dots, \nu_6  \subset W
\]
be the lines containing $p$.   Then $\ell_x$ is the only one of them which is contained in 16 quadrics, and this characterizes $X$ as the component of $F(W)$ containing $[\ell_x]$.  \qed

\begin{corollary}\label{C:pala2}
The rational map
\[
\xymatrix{w:G(5,\bigwedge^2V^\vee)\ar@{-->}[r]&  \U &
 [U] \ar@{|->}[r]& D_4(\phi)}
 \]
 where $\phi: U\otimes\O_{\P(V)} \to \Omega^1_{\P(V)}(2)$ is the evaluation map and
 $\U \subset \P(S^4V^\vee)$ is the locally closed subscheme of the linear system of quartics in $\P(V)$ 
which parametrizes  the Palatini quartics, is birational.
\end{corollary}

\proof 
To each 
 $[U] \in G(5,\bigwedge^2V^\vee)$ we can associate a   linear section
 \[
 X := G(2,V) \cap \P(U)^\perp  \subset \P(\wedge^2V^\vee)
 \]
 This correspondence identifies an open set of  $G(5,\bigwedge^2V^\vee)$ with the family of prime Fano threefolds of genus 8 contained in $G(2,V)$. Therefore the corollary follows from Theorem \ref{C:pala1}. \qed

\begin{remarks}\rm\label{R:pala1}  The  class of $F(W)$ in $A_\ast(G(2,V))$ is easily computed to be (\cite{wF84}, Example 14.7.13 p. 275)
\[
[F(W)] = c_5(S^4Q) = 32c_2(Q)\left[3 c_1(Q)^3+4c_2(Q)c_1(Q)\right]
\]
Consequently the degree of $F(W)$ in $\P(\wedge^2V)$  is
\begin{equation}\label{E:degF}
\begin{array}{ll}
\deg[F(W)] &= c_5(S^4Q)\cdot c_1(Q)^3 \\ \\
&= 96 c_2(Q)c_1^6(Q) + 128 c_2^2(Q)c_1^4(Q)   \\ \\
& = 96 \deg[G(2,5)] + 128 \deg[G(2,4)] \\ \\
& =96\cdot 5 + 128 \cdot 2 = 736
\end{array}
\end{equation}
Using standard notation for Schubert cycles on $G(2,V)$ (see e.g. \cite{wF84}), the class of $F(W)$ can be also written under the form
\[
[F(W)] = 320 (1,3) + 96 (0,4)
\]
where $96= 24 \times 4$ is the number of lines of $F(W)$ meeting a general line of $\P(V)$, and $320$ is the degree of the ruled surface spanned by the lines of $W \cap H$ for a general hyperplane $H$. We have already seen (Prop. \ref{P:lines2}) that $X$ is identified with an irreducible component of $F(W)$, which has degree 14. Its class is 
\[
[X] = 5 (1,3) + 4 (0,4)
\]
where 5 is the degree of the ruled surface spanned by the lines $\ell_x$ which are contained in a general hyperplane $H \subset \P(V)$, namely the degree of the grassmannian of lines of $H$. On the other hand 4  is the degree of $W$,  which is the number of lines parametrized by $X$ and meeting a general line $\ell \subset \P(V)$.

\noindent
Similarly the component $h(Y)$ has class:
\[
[h(Y)] = 10 (1,3) + 4 (0,4)
\]
as one easily computes in a similar way.

\noindent 
Consider a general point $p \in W$ and let
 \[
\ell_x, n_y,
\lambda_1, \dots, \lambda_{16},
\nu_1, \dots, \nu_6  \subset W
\]
be the lines containing $p$ (notation as in the proof  of Prop. \ref{P:pala3}).  

\separation
\noindent
\emph{Claim:} $\lambda_1, \dots, \lambda_{16}$ belong to the same irreducible component of $F(W)$, which we call $X'$, and $\nu_1, \dots, \nu_6$ belong to another irreducible component of $F(W)$ which we call $Y'$.  

\separation
\emph{Proof of the Claim.}  Since the Fano scheme $F_2(X)$ of conics on $X$ is irreducible and nonsingular of dimension two (\cite{enc99}, \S 4.2), the quadrics  $Q_1, \dots, Q_{16}$ belong to an irreducible two-dimensional family of quadrics, and therefore the lines contained in these quadrics are contained in at most two irreducible components of $F(W)$. One of them is $X$, therefore there is only another one, and this is $X'$. 

\noindent
The argument for $Y'$ is similar, relying on the irreducibility of the Fano scheme of lines on the cubic threefold $Y$.

\separation
The components $X,h(Y),X',Y'$ of $F(W)$ are all the irreducible components parametrizing lines going through the general point $p \in W$. There is another irreducible component $D$ of $F(W)$, which
 is supported on the set of $[\lambda] \in G(2,V)$ such that $\lambda$ is contained in a plane $\Pi\subset W$  of the form $p_W(p_x^{-1}(\ell))= \Pi_p$ for some line $\ell \subset X$ and for some $p \in \Gamma(W)$ (the notations are those of Prop. \ref{P:lines1}). Clearly $D$ is irreducible of dimension three and therefore it is an irreducible component of $F(W)$. One easily finds that its class is $[D]= 45 (1,3)$. In particular $[D]$ has degree 90.  
 There are no other components of $F(W)$. Hence:
 \[
 F(W) = X \cup h(Y) \cup X' \cup Y' \cup D
 \]
 
 \end{remarks}

 \section{The main theorem}

  We can summarize all the above  in the following:
  
  \begin{theorem}\label{T:main1}
 By associating to a   Fano threefold    $X =G(2,V) \cap \P(U)^\perp$,  where $[U] \in  G(5,\bigwedge^2V^\vee)$,  the singular curve
 $\Gamma(W)\subset \P(V)$ of its Palatini quartic  $W$, we obtain a generically injective rational map
 \[
 \xymatrix{
 \gamma: G(5,\bigwedge^2V^\vee)\ar@{-->}[r]& \Hilb^{\P(V)}_{25t-25}}
 \]
 whose image is a locally closed subset of dimension 50.
 \end{theorem}
  
   \proof  The   map $\gamma$ is the composition 
   \[
  \xymatrix{G(5,\bigwedge^2V^\vee)\ar@{-->}[r]^-w&  \U \ar[r]^-\sigma&  \Hilb^{\P(V)}_{25t-25}}
  \]
  where $w$ and $\sigma$ are the maps introduced in \S \ref{S:family}. These maps are both generically injective, by Corollaries \ref{C:gentor2} and \ref{C:pala2}. \qed

    \begin{remark}\rm\label{R:main1}
 If $\Gamma = \Gamma(W)$ is the singular curve of a Palatini quartic, then one easily computes that
 \[
 \chi(N_{\Gamma/\P(V)}) = 100
 \]
 and therefore $\dim_{[\Gamma]}\left[\Hilb^{P(V)}_{25t-25}\right] \ge 100$.  It follows that the image of the map $\gamma$  has high codimension in $\H$.  
 \end{remark}

 Theorem \ref{T:main1} can be also stated in the following, more intrinsic, form:
 
 \begin{theorem}\label{T:main2}
 A general prime Fano threefold  $X$ of genus 8 can be reconstructed, up to isomorphism,  from the pair $(\Gamma,L)$, where $\Gamma$ is its Fano curve of lines and $L=\O_\Gamma(1)$ is the theta-characteristic   which gives the natural embedding $\Gamma \subset \P^5$. 
 \end{theorem}
 
 \proof
 Given a general Fano threefold $X=G(2,V) \cap \P(U)^\perp$ as in Theorem \ref{T:main1}, 
 the pair $(\Gamma,L)$ uniquely determines the embedded curve $\Gamma(W)\subset \P(V)$ up to the action of PGL$(V)$, and therefore, by Theorem \ref{T:main1}, also $X$ is uniquely determined by $(\Gamma,L)$ up to the action induced by PGL$(V)$ on $G(2,V)$. But the PGL$(V)$-orbits of Fano threefolds of genus 8 in $G(2,V)$ coincide with their isomorphism classes because they are embedded anticanonically.  \qed

 \begin{remark} \rm Note that for a pair $(\Gamma,L)$ coming from a Fano threefold $X$, the theta-characteristic  $L$ has  $h^0(\Gamma,L)=6$, and it is very likely that there is only one theta-characteristic on $\Gamma$ having this property. In such a case the correspondence 
 $\xymatrix{[X]\ar@{|->}[r]& [\Gamma]}$ would be one-to-one. 
 \end{remark}

   \section*{Appendix}
 
 In this Appendix we give another approach to the proof of Corollary  \ref{C:gentor2}. We keep the same notation and assumptions as before. We assume given a Palatini quartic $W$ such that the curve $\Gamma(W)$ is nonsingular.   We ask the following:
 
 \separation
 
  \emph{Question:}  Is  $W$ the only Palatini quartic whose singular locus is $\Gamma:=\Gamma(W)$?
 
 \begin{lemma}\label{L:gamma1}
  If the map 
 \[\xymatrix{\partial: I_3\otimes  V \ar[r] & H^0(\P^5,\O(2))= S^2V^\vee \\
 F \otimes (\sum a_i{\partial\over \partial X_i}) \ar@{|->}[r]& \sum a_i {\partial F \over \partial X_i}}\]
 is surjective then the curve $\Gamma$ is the singular locus of   a unique Palatini quartic   $W$, in other words  the above question has a positive answer.
 \end{lemma}
 
 \proof
  Assume that $\partial$ is surjective. Then  the 21 partial derivatives
 \[
 {\partial^2W \over\partial X_i\partial X_j}, \qquad 0 \le i \le j \le 5
 \]
 are    linearly independent because they   generate Im$(\partial)$.  
  We have:
 \[
 \dim(\ker(\partial)) = \dim[I_3\otimes  V] - h^0(\P^5,\O(2)) = 15
 \]
 On the other hand, $\ker(\partial)$ contains the  space $D$ generated by the 15 tensors:
 \[
 {\partial W \over \partial X_j}\otimes {\partial\over\partial X_i} - 
 {\partial W \over \partial X_i}\otimes {\partial\over\partial X_j},  
 \qquad 0 \le i < j \le 5
 \]
 Since  they are linearly independent,  $\dim(D)=15$ and $D = \ker(\partial)$.
 
 Assume that there is another Palatini quartic  $W'$ such that $\Gamma = {\rm Sing}(W')$. Then 
 $\ker(\partial)$ is also generated by the tensors:
 \[
 {\partial W' \over \partial X_j}\otimes {\partial\over\partial X_i} - 
 {\partial W' \over \partial X_i}\otimes {\partial\over\partial X_j}, 
 \qquad 0 \le i < j \le 5
 \]
 There are two possibilities:
 \begin{itemize}
 
 \item[(a)] For each $0 \le i < j \le 5$ there is $\alpha_{ij}\in {\bf C}$ such that: 
 \[{\partial W \over \partial X_j}\otimes {\partial\over\partial X_i} - 
 {\partial W \over \partial X_i}\otimes {\partial\over\partial X_j} 
 = \alpha_{ij}\left({\partial W' \over \partial X_j}\otimes {\partial\over\partial X_i} - 
 {\partial W' \over \partial X_i}\otimes {\partial\over\partial X_j}\right)
 \]

 \item[(b)]  For some $i < j$  the tensors 
 \[
 {\partial W \over \partial X_j}\otimes {\partial\over\partial X_i} - 
 {\partial W \over \partial X_i}\otimes {\partial\over\partial X_j}
 \]
    and 
 \[
 {\partial W' \over \partial X_j}\otimes {\partial\over\partial X_i} - 
 {\partial W' \over \partial X_i}\otimes {\partial\over\partial X_j}
 \]
  are linearly independent.
 \end{itemize}
 
 In case (a) we obtain the identities:
 \[
 {\partial (W - \alpha_{ij}W')\over \partial X_j}\otimes {\partial\over\partial X_i} =
 {\partial (W  - \alpha_{ij}W')\over \partial X_i}\otimes {\partial\over\partial X_j} 
 \]
 which imply that  $\alpha_{ij} = \alpha$ are independent of $i,j$ and
 \[
 {\partial W \over \partial X_i} = \alpha {\partial W' \over \partial X_i}
 \]
 for all $i$.  This implies that the two Palatini quartics $W$ and $W'$ are equal because they have the same first polars with respect to every point. 
 
 In case (b) we obtain that the 16 tensors 
 \[
 {\partial W \over \partial X_j}\otimes {\partial\over\partial X_i} - 
 {\partial W \over \partial X_i}\otimes {\partial\over\partial X_j}, \quad
 {\partial W' \over \partial X_j}\otimes {\partial\over\partial X_i} - 
 {\partial W' \over \partial X_i}\otimes {\partial\over\partial X_j},
 \]
 \[
 {\partial W \over \partial X_h}\otimes {\partial\over\partial X_k} - 
 {\partial W \over \partial X_k}\otimes {\partial\over\partial X_h}, \qquad
 0 \le k < h \le 5,  \quad (k,h) \ne (i,j)
 \]
 are linearly independent elements of $\ker(\partial)$, a contradiction. \qed

 \begin{remark}\label{R:catal}\rm  The surjectivity of $\partial$ is equivalent to the property that 
 ${\rm Cat}(W)\ne 0$, i.e. that $W$ is not a zero of the \emph{catalecticant determinant}  or, equivalently, that it is not apolar to any quadric  (see \cite{DK93}, \S 2, for details).
\end{remark}

 Even though we cannot prove that $\partial$ is surjective, we have  experimental confirmation that 
 \emph{$\partial$ is surjective for a general choice of $W$}.  In fact, a simple Macaulay2 program 
 \cite{MAC2} produces a Palatini quartic with random coefficients having linearly independent second partial derivatives. We have reproduced the script of the program below,   written by G. Ottaviani to which we are thankful. This, together with Lemma  \ref{L:gamma1},    gives a positive answer to our question for a general  Palatini quartic $W$. 
 
\begin{verbatim}
S=QQ[a..e,f..k]
R=QQ[a..e]
m=random(R^{6:1},R^{6:0});
m=m-transpose(m);
-- we create random skew symmetric matrix m
n=sub(m,S);
-- n is m in the ring S
pn=matrix{{generators(pfaffians(6,n))}};
-- pn is the cubic pfaffian
p=ideal(matrix(pn|matrix{{f,g,h,i,j,k}}*n));
-- p contains the pfaffian and the conditions to lie in the kernel
ae=sub(ideal(a,b,c,d,e),S)
q=quotient(p,ae);
pal=eliminate({a,b,c,d,e},q);
-- pal is Palatini quartic
T=QQ[f..k]
pala=sub(pal,T);
-- pala is Palatini quartic in the new ring T
jpala=jacobian(pala);
jjpala=jacobian(jpala);
betti trim(ideal(jjpala))
-- output is 21, hence the second derivatives of Palatini quartic are
independent
\end{verbatim}

   \noindent
   \textsc{f. flamini} -
   Dipartimento di Matematica, 
 Universit\`a degli Studi di Roma ``Tor Vergata'',  Viale della Ricerca Scientifica, 1 - 00133 Roma (Italy).
 e-mail: \texttt{flamini@mat.uniroma2.it}

\noindent
\textsc{e. sernesi} - 
 Dipartimento di Matematica,
 Universit\`a Roma Tre, 
 Largo S.L. Murialdo, 1 -
 00146 Roma (Italy).  e-mail: \texttt{sernesi@mat.uniroma3.it}
 
\end{document}